\magnification= \magstep1
%\hsize 6.00 true in
%\vsize 8.0 true in
%\hoffset=.45in
%\voffset=.7in
\baselineskip=14pt
%\nopagenumbers
\pageno=0
\newcount\knum
\knum=1
\newcount\uknum
\uknum=1
\newcount\znum
\znum=1
\newcount\notenumber

\def\Za{\the\knum.\the\uknum.\the\znum \global\advance\znum by 1}
\def\ukneu{\znum=1\global\advance\uknum by 1}
%\headline={\ifodd\pageno\rightheadline \else\leftheadline\fi}
%\def\rightheadline{\vbox{\line{\petit\noindent\hskip -0.7 true cm
%\sl\S 3 Geometrie prinzipaler Einbettungen
%\hfil{\tenrm \folio}}\vskip 0.7cm}}
%\def\leftheadline{\vbox{\line{\petit\noindent\sl\hskip -0.7 true cm
%{\tenrm \folio}\hfil ????}
%\vskip 0.7cm}}
\font\tbfontt=cmbx10 scaled \magstep0
\font\gross=cmbx10 scaled \magstep2
\font\mittel=cmbx10 scaled \magstep1
\font\ittel=cmr10 scaled \magstep1

\font\eightrm=cmr8

\font\eightit=cmti8 scaled \magstep0
\def\sqr#1#2{{\vcenter{\vbox{\hrule height.#2pt\hbox{\vrule width.#2pt
height#1pt \kern#1pt \vrule width.#2pt}\hrule height.#2pt}}}}

\def\cz{{\rm C\hskip-4.8pt\vrule height5.8pt\hskip5.3pt}}
\def\kz{{\rm C\hskip-4.0pt\vrule height4.3pt\hskip5.3pt}}  %Index C-Zeichen

\def\nz{{\rm I\hskip -2pt N}}

  %quaternionisches Zahlenzeichen

\def\lim{\mathop{\rm lim}}

\parindent=0pt
\font\eightrm                = cmr8
\font\eightsl                = cmsl8
\font\eightsy                = cmsy8
\font\eightit                = cmti8

\font\eighti                 = cmmi8
\font\eightbf                = cmbx8
\def\petit{\def\rm{\fam0\eightrm}
\textfont0=\eightrm %\scriptfont0=\sixrm \scriptscriptfont0=\fiverm
 \textfont1=\eighti %\scriptfont1=\sixi \scriptscriptfont1=\fivei
 \textfont2=\eightsy %\scriptfont2=\sixsy \scriptscriptfont2=\fivesy
 \def\it{\fam\itfam\eightit}
 \textfont\itfam=\eightit
 \def\sl{\fam\slfam\eightsl}
 \textfont\slfam=\eightsl
 \def\bf{\fam\bffam\eightbf}
 \textfont\bffam=\eightbf %\scriptfont\bffam=\sixbf
 %\scriptscriptfont\bffam=\fivebf
 \normalbaselineskip=9pt
 \setbox\strutbox=\hbox{\vrule height7pt depth2pt width0pt}
 \normalbaselines\rm}
\newdimen\refindent
\def\begref{\vskip1cm\bgroup\petit
\setbox0=\hbox{[Bi,Sc,So]o }\refindent=\wd0
\let\sl=\rm\let\INS=N}

\def\ref#1{\filbreak\if N\INS\let\INS=Y\vbox{\noindent\tbfontt
References\vskip1cm}\fi\hangindent\refindent
\hangafter=1\noindent\hbox to\refindent{#1\hfil}\ignorespaces}

\long\def\fussnote#1#2{{\baselineskip=9pt
\setbox\strutbox=\hbox{\vrule height 7pt depth 2pt width 0pt}%
\petit\noindent\footnote{\noindent #1}{#2}}}
{\nopagenumbers
%\rightline{CPT-95/P.????}
%\rightline{Mannheimer Manuskripte 181}
\rightline{hep-th/9506???}
%\footline={\hfil}
\rightline{June 1995}
\vskip 4cm
\phantom{prelim.version}
\centerline{\gross A Note on the Wodzicki Residue}
\smallskip
\phantom{\rm (preliminary version)}
\vskip 1cm

\centerline{\ittel Thomas Ackermann\footnote{$^*$}{\eightrm
e-mail: ackerm@euler.math.uni-mannheim.de} }
\vskip 0.5cm
\centerline{Wasserwerkstr. 37, D-68309 Mannheim, F.R.G.}
\vskip 1.5cm
\centerline{\vbox{\hsize=5.0 true in \petit\noindent
\bf Abstract.\rm\ In this note we explain the relationship
of the Wodzicki residue of (certain powers of) an elliptic
differential operator $P$\ acting on sections of a complex vector
bundle $E$\ over a closed compact manifold $M$\
and the asymptotic expansion
of the trace of the corresponding heat operator
$e^{-tP}$. In the special case of a generalized
laplacian $\triangle$\ and ${{\rm dim}\;M > 2}$, we
thereby obtain a simple proof of the fact already shown in
[KW], that the Wodzicki residue $res(\triangle^{-{n\over 2}+1}
)$\ is the integral of the second coefficient of the heat
kernel expansion of $\triangle$\ up to a proportional
factor. }}
\par
1991 MSC: \it 47F05, 53B50, 58G11\hfill \vfill
\break}\rm
\advance\hsize by -0.5 true in
\advance\vsize by -0.9true in
\advance\hoffset by 0.45 true in
\advance\voffset by 0.7 true in
{\mittel 1. Introduction}
%\ukneu
\vskip 0.7cm
In [KW] it was shown that the Wodzicki
residue $res(\triangle^{-{n\over 2}+1})$\
of a generalized laplacian
$\triangle$\ on a complex vectorbundle $E$\ over a closed compact
manifold $M$, ${\rm dim}\;M=n>2$, is
the integral of the second coefficient of the heat
kernel expansion of $\triangle$\ up to a proportional
factor. This observation reflects a more
general property which we explain
in this note, namly that the Wodzicki residue  
of (certain powers of) an elliptic
differential operator $P$\ on $E$\
can be related with the coefficients in the asymptotic expansion
of the trace of the corresponding heat operator
$e^{-tP}$.  
Even if this might be well-known by mathematicians working in
the field - our simple proof relies on two results made by
Wodzicki [W] and Gilkey [Gi] - we think it worth to restate it
because of the growing importance of the Wodzicki residue
for non commutative geometry, for example to incorporate
gravity in the Connes-Lott model (cf. [CFF]). By the way
we correct an error term in [KW].
\vskip 1.0cm
{\mittel 2. Computing ${res(P^{-({n-k\over d}})}$
 }
\vskip 0.7cm
\ukneu
Let $E$\ be a finite dimensional complex vector bundle over
a closed compact manifold $M$\ of dimension $n$. Recall that the
non-commutative residue of
a pseudo-differential operator $P\in \Psi {\rm DO}(E)$\
can be defined by
$${res(P):={(2\pi)^{-n}}\;\int_{S^*M}\;
tr\bigl(\sigma^P_{-n}(x,\xi)\bigr)\;dx d\xi,}\eqno(2.1)$$
where $S^*M\subset T^*M$\ denotes the co-sphere bundle on $M$\ and
$\sigma^P_{-n}$\ is the component of order $-n$\ of the complete
symbol $\sigma^P:=\sum_i\;\sigma^P_i$\ of $P$, cf [W]. 
In his thesis, Wodzicki has shown that the linear
functional $res\colon
\Psi {\rm DO}(E,F)\rightarrow \cz$\ is in fact
the unique trace (up to multiplication by constants) on
the algebra of pseudo-differential operators $\Psi{\rm DO}(E)$.
\smallskip
Now let $P\in \Psi DO(E)$\ be elliptic with ${\rm ord}\;P
=d>0$. It is
well-known (cf. [Gi]) that its zeta function $\zeta(P,s)$\ is a
holomorphic function on the half-plane ${{\rm Re}\;s > n/d}$\ which
continues analytically to a meromorphic function
on $\kz$\ with
simple poles
at \vfil\break
$\{\;{(n-k)\over d}\;\vert k\in \nz\setminus \{ n\}\;\}$.
For $n-k > 0$\ with $k\in \nz$\ we have the following identities:
$$\eqalignno{res(P^{-({n-k\over d})}) &= d\cdot Res_{s={n-k\over d}}
\zeta(P,s) &(2.2)\cr
Res_{s={n-k\over d}}\zeta(P,s)&= a_k(P)\cdot \Gamma(
\hskip -0.6cm\hbox{\petit ${n-k\over d}$})^{-1}.
&(2.3) \cr}$$
Here $\Gamma$\ is the gamma function and $a_k(P)$\ denotes the
the coefficient of $t^{{k-n\over d}}$\ in the
asymptotic expansion of $Tr_{L^2}e^{-tP}$. The first equality was
shown in [W], whereas (2.3) is a consequence of the Mellin transform
which relates the zeta function and the heat equation and was already
proven by Gilkey in [Gi]. Consequently we obtain
$${a_k(P)= d^{-1}\cdot\Gamma(
\hskip -0.6cm\hbox{\petit ${n-k\over d}$})\cdot 
res(P^{-({n-k\over d})}).}\eqno(2.4)$$
Now suppose $P=\triangle$\ is a generalized laplacian and $k=2$.
Then $a_2(\triangle)=(4\pi)^{-{n\over 2}}\phi_2(\triangle)$\ where
$\phi_2(\triangle)$\ denotes the integral over the diagonal part
of the second coefficient of the heat kernel expansion of
$\triangle$. By
using the well-known identity $\Gamma(z+1)=z\Gamma(z)$\ it is
easily verified that 
$${(n-2)\phi_2(\triangle)=(4\pi)^{n\over 2}\cdot
\Gamma(
\hskip -0.6cm\hbox{\petit ${n\over 2}$})\cdot 
res(\triangle^{-{n\over 2}+1})}\eqno(2.5)$$
for $n>2$. Note, that in the above mentioned
reference [KW] the Wodzicki
residue was defined
by
${\widetilde{res}(P)\!:=\!{\Gamma({n\over 2})\over 2
\pi^{n\over 2}}\int_{S^*M}\;
tr\bigl(\sigma^P_{-n}(x,\xi)\bigr) dx d\xi}$,
using $vol(S^{n-1})^{-1}={\Gamma({n\over 2})\over 2
\pi^{n\over 2}}$\ as a nomalizing factor.
Thus, equation (2.5)
can be equivalently expressed as
$\widetilde{res}(\triangle^{-{n\over 2}+1})
={n-2\over 2}\phi_2(\triangle)$. Therefore their proportional
factor in equation (4.10) is not correct.
\smallskip
We conclude in exploiting the result we have in hand
in the case of a Clifford module
$E={\cal E}$\ over a four-dimensional
closed compact Riemannian manifold $M$\ and
$\triangle$\ being the square
of a Dirac operator
$D$\ defined by a Clifford connection.
Then
equation (2.5) together with the
well-known observation $\phi_2(D^2) \sim \int_M *r_M$,
where $*$\ is the Hodge star corresponding to the Riemannian
metric and $r_M$\ denotes the scalar curvature of $M$,
yields $res(D^{-2})\sim \int_M\; *r_M$. Thus, by the above
derivation of (2.5),  
the fact - first announced by Connes [C] and shown
by Kastler [K] using the symbol calculus -
that the non-commutative residue of the inverse square of
the Dirac operator is proportional to the Einstein-Hilbert action of
general relativity, is almost obvious.
\vskip 0.5cm
{\petit\bf Note:\rm\ After completion of this work we have been told
that this result was independently recognized by
M. Walze [Wa].}\vfil\break
\begref
\ref{[C]} A.Connes, \sl Non-commutative geometry and physics\rm ,
IHES preprint (1993)
\ref{[CFF]} A.H.Chamseddine, G.Felder, J.Fr\"ohlich, \sl
Gravity in non-commutative geometry\rm , Com. Math. Phys {\bf 155} (1993),
205-217
\ref{[Gi]} P.G.Gilkey, \sl Invariance theory, the heat equation,
and the Atiyah-Singer index theorem, \rm Publish or
Perish (1984)
\ref{[K]} D.Kastler, \sl The Dirac operator and gravitation\rm ,
Com. Math. Phys. \bf 166\rm (1995)
\ref{[KW]} W.Kalau, M.Walze, \sl Gravity, non-commutative geometry
and the Wodzicki residue\rm , to appear in Journ. of
Geometry and Physics
\ref{[W]} M.Wodzicki, \sl Non-commutative residue I\rm , LNM {\bf 1289}
(1987), 320-399
\ref{[Wa]} M.Walze, \sl Nicht-kommutative Geometrie und
Gravitation, thesis, Universit\"at Mainz, to appear 
\bye
Moreover, in [KW] - also
using the symbol calculus - equation (3.11) was proven\fussnote{${^{(
1)}}$}{Using $vol(S^{n-1})^{-1}={\Gamma({n\over 2})\over 2
\pi^{n\over 2}}$\ as a nomalizing factor in [KW] the non-commutative
residue was defined as
${\widetilde{res}(P):={\Gamma({n\over 2})\over 2\pi^{n\over 2}}\;\int_{S^*M}\;
tr\bigl(\sigma^P_{-n}(x,\xi)\bigr)\;dx d\xi}$. With this definition
(3.11) can be restated as $\widetilde{res}(\triangle^{-{n\over 2}+1})
={n-2\over 2}\phi_2(\triangle)$. So the proportionality constant
in [KW], equation (4.10) is false.}.
\bye